\documentclass[]{article}
\usepackage[utf8]{inputenc}
\usepackage{amssymb}
\usepackage{amsmath}
\usepackage{amsthm}
\usepackage{epsfig}
\usepackage{graphics}
\usepackage{latexsym}
\usepackage{rotating}
\usepackage{fancyhdr}
\usepackage{textcomp}
\usepackage{tikz}
\usepackage{float}
\usepackage{indentfirst}

\numberwithin{equation}{section}
\newtheorem{teo}{Theorem}[section]
\newtheorem{tnm}{Definition}[section]
\newtheorem{lmm}{Lemma}[section]
\newtheorem{uyr}{Remark}[section]
\newtheorem{snc}{Result}[section]

\title{{\textbf{Time-fractional quenching problem: Blow-up of $D_{t}^{\alpha}u$ at the quenching point } }}
\author{Nurdan Kar \\ \scriptsize{Department of Mathematics, Ankara University Ankara, 06100, Turkiye} \\ 
\scriptsize{Department of Mathematics, University of Maryland, Maryland, 20742, USA} \\\scriptsize{E-mail: kar@ankara.edu.tr, nkar1@umd.edu}}
\date{}

\begin{document}
\maketitle
\begin{abstract}
In this paper, we establish the occurrence of blow-up in the time-fractional term at the quenching point. By demonstrating that the quenching points are contained within a compact subset of the designated spatial interval, we employ this finding to prove the blow-up of the Caputo fractional time-derivative at the quenching point.\\

\end{abstract}

\section{Introduction}
In recent decades, there has been a growing interest in fractional differential equations (FDEs) due to their potential applications in various fields of science and engineering, as well as their theoretical significance \cite{das}, \cite{kilbas}, \cite{podlubny}, \cite{oldham}. Within the realm of FDEs, considerable attention is devoted to fractional diffusion equations, owing to their paramount importance in capturing the essence of physical phenomena characterized by memory effects \cite{caputo}. Caputo introduced a distinct memory mechanism that is employed as a time derivative in this study \cite{caputo2}, \cite{caputo}.

Although the comprehension of linear diffusion equations, often used to represent ideal diffusive phenomena, is relatively straightforward, the behavior of solutions becomes considerably more sophisticated when dealing with nonlinear (conventional or fractional) diffusion equations. In 1975, Kawarada proposed an interesting one-dimensional initial-boundary value problem denoted as 
\begin{equation} \label{main}
		\begin{aligned}
			\begin{cases}
				\dfrac{\partial}{\partial t}u(x,t) -\dfrac{\partial^{2}}{\partial x^2}u(x,t)  = \dfrac{1}{1-u} &\quad \text{in} \, (0,L)\times (0,T),  \\[2mm]
				u(0, t)=u(L, t)=0 &\quad \text{in}\, (0,T), \\[2mm]
				u(x, 0)=0  &\quad \text{on}\, [0,L], 
			\end{cases}
		\end{aligned}
	\end{equation}
where $L$ is a positive real number indicating the length of spatial domain \cite{kawarada}. Subsequently, he achieved remarkable outcomes as follows:
\begin{snc}
If the solution $u(x,t)$ reaches $1$ within a finite time, then $u_t(x,t)$ is unbounded. That is,
\begin{equation}\nonumber
\lim _{t \rightarrow T^{-}} \sup \left\{\mid u_{t}(x, t) \mid: x \in[0, L]\right\} \rightarrow \infty.
\end{equation}
Here, the necessary condition for quenching is given by
\begin{equation}\nonumber
		\lim _{t \rightarrow T^{-}} \max \{\mid u(x, t)\mid : x \in[0, L]\} \rightarrow 1^{-}.
\end{equation}
The value $T$ represents the quenching time.
\end{snc}
\begin{snc}
If $L > 2\sqrt{2}$ in the problem (\ref{main}), $u(L/2,t)$ reaches one in finite time.
\end{snc}

It is easy to see that the problem (\ref{main}) represents a typical one-dimensional diffusion model with a singular source term. In contrast to our model in this study, the temporal and spatial derivatives in problem (\ref{main})  are classic first and second order derivatives, respectively. However, it can be readily stated that the growing interest in fractional analysis has extended to encompass the quenching problem. In \cite{xu1} and \cite{xu2018}, the fractional quenching time has been estimated, and some numerical simulations have been presented for the time-fractional quenching problem with the Caputo derivative.  In \cite{xu2019}, it has been focused on a two-dimensional time-fractional quenching problem, and a novel method has been proposed to enhance the accuracy of estimating quenching moments. There are also studies in the literature using the fractional derivative acting on the spatial variable, see \cite{xu2017a} and references therein.

The investigation of the relation between quenching and blow-up phenomena has been another focal point within the framework of the (\ref{main}) problem type \cite{caffarelli}, \cite{deng}, \cite{friedman}. In \cite{ke}, the blow-up of classic time derivative at the quenching point in degenerate parabolic equations has been examined. Besides, the exploration of classic time-derivatives blowing up at the quenching point has been studied in the context of a reaction-diffusion system with logarithmic singularity, as studied in \cite{mu}. However, although numerical studies have shown the occurrence of blow-up at the quenching point in fractional quenching problems \cite{xu1}, \cite{xu2018}, \cite{xu2019} a theoretical proof for this phenomenon is currently unavailable. In this study, we provide a proof for this occurence.

The importance of the maximum principle for parabolic differential equations is well recognized. In particular, the maximum principle for the generalized time-fractional diffusion equation (TFDE) with the Caputo–Dzherbashyan fractional derivative has been established in \cite{luchko2009} through an extremum principle, over an open bounded domain $G \times (0, T)$ where $G \subset \mathbb{R}^n$. It is worth noting that there exists a maximum principle for a suitably defined weak solution in the fractional Sobolev spaces, although not applicable to the strong solution \cite{luchko2017}. We leverage this principle effectively to obtain the results presented in the subsequent sections of this paper.

In this study, inspired by Kawarada's research, we investigate the following time-fractional problem
\begin{equation}\label{m}
	\begin{aligned}
		\begin{cases}
			D_t^{\alpha}u(x,t) - \dfrac{\partial^{2}}{\partial x^2}u(x,t) = f(u)  & \quad \text{in} \, (0,L)\times(0,T),  \\[2mm]
			u_x(0,t)=0,  u(L,t)=0  &\quad \text{in} \, \left[0,T\right), \\[2mm]
			u(x,0) =u_0(x)  & \quad \text{on} \, [0,L],
		\end{cases}
	\end{aligned}
\end{equation}
where $0<\alpha<1$, the operator $D_t^{\alpha}$ is the Caputo fractional derivative of order $\alpha$, and $f(u)$ is a nonlinear source term with singularity on $u$. We broaden the scope of the classical quenching equation by incorporating a time-fractional formulation utilizing the Caputo derivative. Furthermore, we impose a derivative on the left boundary condition in problem (\ref{m}). It is clear that the main theoretical results demonstrated here are also applicable to TFDE with different boundary conditions.

The remaining sections of this paper are structured as follows: Section 2 presents a concise overview of the definition and several properties of fractional derivatives. Chapter 3 focuses on presenting and proving lemmas and theorems that play a crucial role in supporting the main results of this study's problem. Section 4 presents the main findings regarding the relationship between blow-up and quenching phenomena in the TFDE.

\section{Mathematical preliminaries}
This section presents an exposition of fundamental definitions and properties pertaining to fractional integrals and derivatives. Comprehensive discussions on these topics can be found in various references \cite{das},  \cite{kilbas}, \cite{podlubny}, \cite{samkokilbas}, \cite{oldham}. Throughout this paper, let $\Gamma(\cdot)$ be Euler's gamma function. 

\begin{tnm} \cite[page 69]{kilbas}
Let $[a,b]$ be a finite interval and $\alpha >0$. Then the $\alpha$th order left-sided and right-sided Riemann–Liouville fractional integrals of $u\in L^{1}[a,b]$ are respectively defined as
\begin{equation}\nonumber
	\begin{aligned}
		&\left({ }_{a} I_{x}^{\alpha}\right)(u(x))=\frac{1}{\Gamma(\alpha)} \int_{a}^{x}(x-s)^{\alpha-1} u(s) d s, \quad x>a, \\[2mm]
		&\left({ }_{x} I_{b}^{\alpha}\right)(u(x))=\frac{1}{\Gamma(\alpha)} \int_{x}^{b}(s-x)^{\alpha-1} u(s) d s, \quad x<b.
	\end{aligned}
\end{equation}
\end{tnm}
\begin{tnm}\cite[page 35]{samkokilbas}
Given a function $u(x)$ defined on a finite interval $[a,b]$, the left and right Riemann-Liouville fractional derivatives of order $\alpha$, $0<\alpha<1$ are respectively denoted by 
\begin{equation}\nonumber
	\begin{gathered}
		\left({ }_{a} D_{x}^{\alpha}\right)(u(x))=\frac{1}{\Gamma(1-\alpha)} \frac{d}{d x} \int_{a}^{x} \frac{u(s) d s}{(x-s)^{\alpha}}, \quad x>a, \\[2mm]
		\left({ }_{x} D_{b}^{\alpha}\right)(u(x))=-\frac{1}{\Gamma(1-\alpha)} \frac{d}{d x} \int_{x}^{b} \frac{u(s) d s}{(s-x)^{\alpha}}, \quad x<b.
	\end{gathered}
\end{equation} 
\end{tnm}
\begin{tnm}\cite[page 92]{kilbas}
Let $0<\alpha<1$  and $u(x) \in AC[a,b]$. Then the left-sided and right-sided Caputo fractional derivatives of $u$ exist almost everywhere and are respectively defined as
\begin{equation}\nonumber
	\begin{aligned}
		&\left({ }_{a}^{c} D_{x}^{\alpha}\right)(u(x))=\left({ }_{a} I_{x}^{1-\alpha}\right) \frac{d u(x)}{d x}, \quad x>a, \\[2mm]
		&\left({ }_{x}^{c} D_{b}^{\alpha}\right)(u(x))=-\left({ }_{x} I_{b}^{1-\alpha}\right) \frac{d u(x)}{d x}, \quad x<b.
	\end{aligned}
\end{equation}
\end{tnm}
The set $AC[a,b]$ denotes the collection of all functions that are absolutely continuous in domain $[a,b]$.

Similarly, the definition of a partial fractional derivative for a given multiple-variable function can be expressed as follows. 
\begin{tnm} \cite[page 358]{kilbas}
The Caputo fractional derivative of a function $u(x,t)$, with $0<\alpha<1$ denoting the order of the derivative, is defined as
\begin{equation}\nonumber
\left({ }D_{t}^{\alpha}\right)(u(x,t))=\frac{1}{\Gamma(1-\alpha)} \int_{0}^{t}(t-s)^{-\alpha} \frac{\partial u(x,s)}{\partial s} d s
\end{equation}
where $0<x<L$ and $t>0$.
\end{tnm}
If $u(x,t)$  is continuously differentiable on  $t$, then $\left({ }D_{t}^{\alpha}\right)(u(x,t))\rightarrow \dfrac{\partial u(x,t)}{\partial t}$  as $\alpha \rightarrow 1$.

\section{The occurrence of quenching}

In this section, we present certain results derived from \cite{xu1}, which are conveniently modified and adapted to our specific problem for the sake of continuity in this paper.

By applying the method of separation of variables \cite{gejji}, \cite{jafari}, the linear part of equation (\ref{m}) can be divided into two ordinary differential equations.  Assume that
\begin{equation}\nonumber
u(x,t)=\Phi(x)\Psi(t).
\end{equation}
Then we have
\begin{equation}\label{özd}
	D_{t}^{\alpha} \Psi(t)+\lambda^{2} \Psi(t)=0, \quad \Psi(0)=u_0(x)
\end{equation}
and
\begin{equation}\label{özd2}
	\Phi^{\prime \prime}(x)+\lambda^2\Phi(x)=0, \quad \Phi^{\prime}(0)=\Phi(L)=0 
\end{equation}  
where $\lambda >0$ is the eigenvalue. Through direct calculation, it can be determined that all eigenvalues are given by the form $\lambda_{n}= (2n-1)\pi / 2L,  n=1,2,... \,. $
Since we are interested in the case of $n=1$,  we  obtain $\lambda_{1}=\pi/2L$ and its associated eigenfunction is $\Phi_1(x)=A\cos(\pi x/2L)$ where  $A$ is an arbitrary constant. There is a particular $A$, such that 
\begin{equation}\label{öz6}
	\int_{0}^{L}\Phi_1(x)dx=1, \quad \Phi_1(x)>0 \, \,\, \text{in} \, (0,L).
\end{equation}
Additionally, the following assumptions regarding (\ref{m})  are necessary for the sake of convenience in the subsequent discussion.

\begin{enumerate}
	\item $f(u)$ is locally Lipschitz.
	\item $f(0)>0$ and $f\in C^1((0,T))$, $f>0$, $f^{'}\geq 0$, $f^{''}\geq 0$.
	\item $\lim _{u \rightarrow u^{{*}^{-}}} f(u)=+\infty$, $u^{*}$ is an isolate singular point.
	\item $f(u)$ satisfies growth condition: $f(u) \geq c_{1}+c_{2}u$, $u\in \left[0,u^{*}\right)$, $c_{1}, c_{2} \in \mathbb{R}^{+}$.
	\item  Considering the classical solution $u(x,t)$ before the quenching moment, it holds that $0 \leq u_{0} < u(x,t) < u^{*}$ for $0 < t < T$ and $0 < x < L$.
	\item  $\int_{0}^{L} \Phi_{1}(x) u_{0}^{ 2}(x) \mathrm{d} x<\infty$ with $u_{0}^{\prime}(0)=0=u_{0}(L)$.
\end{enumerate}

\begin{lmm}\label{lmm2y}\cite{xu1}
Let $\Phi_1(x)$ be the first eigenfunction (\ref{özd2}) and $u \in \mathcal{C}^{2}(0, L) \times A C[0, T]$ then
\begin{equation}\nonumber
\int_{0}^{L} D_{t}^{\alpha} u(x, t) \Phi_{1}(x) d x=D_{t}^{\alpha} \int_{0}^{L} u(x, t) \Phi_{1}(x) d x.
\end{equation}
\end{lmm} 
The subsequent theorem is obtained from \cite{xu1} with certain modifications.
\begin{teo}\label{teosnm2}
Let  $u(x,t)$ be the classic solution of problem (\ref{m}). If the source term $f(u)$ satisfies assumptions (1)-(4) and $\lambda_1^{2}< c_{2}$, then the solution $u(x,t)$ shall quench within a finite time.
\end{teo}
	
	\noindent \textbf{Proof.} Assume that $\left[0, T\right)$ is the largest interval such that $u(x,t)$ continuously well defined. Thus, according to the maximum principle, assumptions (1)-(3), we obtain  
\begin{equation}\nonumber	
	u(x,t)>0, 0<x<L, 0<t<T. 
\end{equation}	
Due to the isolated singularity at point $u^{*}$, we have $\lim _{t \rightarrow T} \sup u(x, t)=u^{*}$.
	
By multiplying both sides of the equation in (\ref{m}) by  $\Phi_{1}(x)$,  and integrating from $0$ to $L$, we have
	\begin{equation}\nonumber
		\begin{aligned}
			\int_{0}^{L} D_{t}^{\alpha} u(x,t) \Phi_{1}(x) d x &=\int_{0}^{L} \dfrac{\partial^{2}}{\partial x^2}u(x,t)  \Phi_{1}(x) d x+\int_{0}^{L} f(u)  \Phi_{1}(x) d x \\[2mm]
			&=\int_{0}^{L} u(x,t)  \dfrac{\partial^{2}}{\partial x^2} \Phi_{1}(x) d x+\int_{0}^{L} f(u)  \Phi_{1}(x) d x \\[2mm]
			&=-\lambda_{1}^{2} \int_{0}^{L} u(x,t) \Phi_{1}(x) d x+\int_{0}^{L} f(u) \Phi_{1}(x) d x 
		\end{aligned}
	\end{equation}
Based on the consideration of Lemma \ref{lmm2y} and assumption (4), we can conclude that
	\begin{equation}\nonumber
		D_{t}^{\alpha}  \int_{0}^{L} u(x,t) \Phi_{1}(x) d x 	\geq-\lambda_{1}^{2} \int_{0}^{L} u(x,t) \Phi_{1}(x) d x+\int_{0}^{L} (c_1+ c_2u(x,t)) \Phi_{1}(x) d x. 
	\end{equation}
By denoting $ y(t)=\int_{0}^{L} u(x, t) \Phi_{1}(x) d x$ and considering (\ref{öz6}), it follows that
	\begin{equation}\label{5}
		D_{t}^{\alpha} y(t) \geq y(t)(c_{2}-\lambda_{1}^{2})+c_{1}
	\end{equation}
with
	\begin{equation}\label{6}
		y(0)=\int_{0}^{L} u(x, 0) \Phi_{1}(x) dx=m_{0} <u^{*}<+\infty.
	\end{equation}
Equations (\ref{5}) and (\ref{6}) represent a fractional initial value problem. 
According to \cite{luchko19991}, \cite{luchko19992}, the solution to (\ref{5}) and (\ref{6}) can be formally expressed as
\begin{equation}\label{çözüm}
		y(t) \geq \int_{0}^{t}\tau^{\alpha-1} E_{\alpha, \alpha}\left[(c_{2}-\lambda_{1}^{2})\tau^{\alpha}\right] d\tau+y(0) E_{\alpha}\left[(c_{2}-\lambda_{1}^{2})t^{\alpha}\right].
	\end{equation}
Taking the limit of $t \rightarrow T$ on both sides  of (\ref{çözüm}) shows that
	\begin{equation}\nonumber
		u^{*} \geq \int_{0}^{T}\tau^{\alpha-1} E_{\alpha, \alpha}\left[(c_{2}-\lambda_{1}^{2})\tau^{\alpha}\right] d\tau+m_{0} E_{\alpha}\left[(c_{2}-\lambda_{1}^{2})T^{\alpha}\right]. 
	\end{equation}
Since $u^{*}$ is finite, it implies that the time $T$ must also be finite. $\qedsymbol$

\section{Main results}	
In this section, we determine the location of the quenching points and prove that $D_t^{\alpha}u(x,t) \rightarrow \infty$  as $u(x,t) \rightarrow u^{{*}^{-}}$, where $u^{*}$  is a singular point. 

\begin{tnm}\label{czmler2}
 $\overline{u}(x,t)$  is called an upper solution of problem (\ref{m}) if  $\overline{u}(x,t)$ satisfies $\overline{u}(0, t)\geq 0$ and the following condition:
\begin{equation}\nonumber
\begin{aligned}
\begin{cases}
D_t^{\alpha}\overline{u}(x,t)-\dfrac{\partial^{2}}{\partial x^2}\overline{u}(x,t)\geq f(\overline{u}) & \quad \mathrm{in} \,(0, L) \times(0, T), \\[2mm]
\overline{u}_x(0, t)\geq 0, \quad \overline{u}(L, t)\geq0 & \quad \mathrm{in} \, (0, T), \\[2mm]
\overline{u}(x, 0)\geq u_{0}(x) & \quad \mathrm{on} \, [0, L].
\end{cases}
\end{aligned}
\end{equation}
Similarly $\underline{u}(x,t)$ is called a lower solution of problem (\ref{m}) if $\underline{u}(x,t)$ satisfies $\underline{u}(0, t)= 0$ and the following condition:
\begin{equation}\nonumber
\begin{aligned}
\begin{cases}
D_t^{\alpha}\underline{u}(x,t)-\dfrac{\partial^{2}}{\partial x^2}\underline{u}(x,t)\leq f(\underline{u}) & \quad \mathrm{in} \, (0, L) \times(0, T), \\[2mm]
\underline{u}_x(0, t)= 0, \quad \underline{u}(L, t)=0 & \quad \mathrm{in} \, (0, T), \\[2mm]
\underline{u}(x, 0)\leq u_{0}(x) & \quad \mathrm{on} \, [0, L].
\end{cases}
\end{aligned}
\end{equation}
\end{tnm}	
	
If $u_0(x)$ is a lower solution, we have the following results.	

\begin{lmm}\label{lmms0}
	If  $-u_{0}^{''}(x) \leq f\left(u_{0}(x)\right)$ in $(0, L)$, then the solution $u(x,t)$ of problem (\ref{m}) increases in $t$  in $(0,L)\times (0,T)$, and $D_t^{\alpha}u>0$ in $(0,L)\times (0,T)$. Furthermore, if $-u_{0}^{''}(x) < f\left(u_{0}(x)\right)$ somewhere in $(0,L)$, then $D_t^{\alpha}u>0$ in $(0,L)\times (0,T)$.
\end{lmm}

\noindent \textbf{Proof.} Since $-u_{0}^{''}(x) \leq f\left(u_{0}(x)\right)$ in $(0, L)$  and $u_{0}^{\prime}(0)=0=u_{0}(L)$, it follows that $u_0(x)$ is a lower solution of problem (\ref{m}). The strong maximum principle implies that
$$ 
u(x,t)\geq u_0(x) \quad \text{in} \, (0,L)\times (0,T).
$$

For any $(x_1,t_1)\in (0,L)\times (0,T)$, there exists a $t_2\in (0,T)$ such that $(x_1,t_1)\in(0,L)\times (0,t_2)$. Let $v(x,t)=u(x,t+h)$ with $0<h<T-t_2$, and $w(x,t)=v(x,t)-u(x,t)$. Then, taking into account the mean value theorem, 
\begin{equation}\nonumber
	\begin{array}{ll}
		\begin{cases}
			D_t^{\alpha}w(x,t)-\dfrac{\partial^{2}}{\partial x^2}w(x,t)=f^{\prime}(\xi)w & \quad \text{in} \, (0, L) \times(0, t_2), \\[2mm]
			w_x(0, t)=w(L, t)=0 & \quad \mathrm{in} \, (0, t_2), \\[2mm]
			w(x, 0)=u(x,h)-u_{0}(x) & \quad \mathrm{on} \, [0, L],
		\end{cases}
	\end{array}
\end{equation}
where $\xi$ is between  $u$ and $v$.
Since $f^{\prime}(\xi)$ is bounded and $u(x,h)-u_{0}(x)\geq 0$, the strong maximum principle implies that $w\geq 0$ in $(0, L) \times(0, t_2)$. Hence, 
$$
v\left(x_{1}, t_{1}\right)-u\left(x_{1}, t_{1}\right)=u\left(x_{1}, t_{1}+h\right)-u\left(x_{1}, t_{1}\right) \geq 0.
$$
Since $h$ is arbitrary small and  $(x_1,t_1)$ is an arbitrary point in $(0,L)\times(0, T)$, it follows that $u(x,t)$ increases in $t$ in $(0, L) \times(0, T)$. Thus, $D_t^{\alpha}u(x,t)\geq 0$ in $(0,L)\times(0, T)$.

If  $-u_{0}^{''}(x) < f\left(u_{0}(x)\right)$ somewhere in $(0, L)$, then from the above proof we have $u(x,h)>u_{0}(x)$ in $(0,L)\times(0, T)$, and hence $w(x,t)>0$ in $(0, L) \times(0, t_2)$. Therefore, $u(x,t)$ strictly increases in $t$ in $(0,L)\times(0, T)$.
For any $(\xi, \eta) \in (0,L)\times(0, T)$, there exists a subset $[x_0,x_2]\times[t_0,t_2]$ of $(0,L)\times(0, T)$ such that  $(\xi, \eta) \in (x_0,x_2)\times(t_0,t_2)$. Since
\begin{equation}\nonumber
\begin{array}{ll}
    D_t^{\alpha}(D_t^{\alpha}u(x,t))-\dfrac{\partial^{2}}{\partial x^2}(D_t^{\alpha}u(x,t)) = f^{\prime}(u)D_t^{\alpha}u(x,t) & \text{in } (x_0,x_2)\times(t_0,t_2), \\[2mm]
    D_t^{\alpha}u(x,t) \geq 0 \quad \text{on } [x_0,x_2]\times[t_0,t_2]&
\end{array}
\end{equation}
and $f^{\prime}(u)$ is bounded on $[x_0,x_2]\times[t_0,t_2]$, the maximum principle implies that either $D_t^{\alpha}u(x,t)>0$ or $D_t^{\alpha}u(x,t)\equiv 0$ in $(x_0,x_2)\times(t_0,t_2)$. Since $D_t^{\alpha}u(x,t)\equiv 0$ contradicts the fact that $u(x,t)$ is strictly increasing in $t$, it follows that $D_t^{\alpha}u(x,t)>0$ in $(0,L)\times (0,T)$. Because $(\xi,\eta)$ is arbitrary in $(0,L)\times (0,T)$, we have $D_t^{\alpha}u(x,t)>0$ in $(0,L)\times (0,T)$. $\qedsymbol$

\begin{lmm}\label{lmms1}
If  $-u_{0}^{''}(x) \leq f\left(u_{0}(x)\right)$ in $(0, L)$ and the strict inequality holds somewhere in $(0, L)$, then for any subset $\left[x_{0}, x_{1}\right] \times\left[t_{0}, T\right) \subset(0, L) \times(0, T)$, there exist a constant  
$c_1>0$ such that $D_t^{\alpha}u(x,t) \geq c_1>0$ in $\left[x_{0}, x_{1}\right] \times\left[t_{0}, T\right)$ .
\end{lmm}	

\noindent \textbf{Proof.} Let us consider the problem
\begin{equation}\nonumber
	\begin{aligned}
		\begin{cases}
			D_t^{\alpha}v(x,t) -\dfrac{\partial^{2}}{\partial x^2}v(x,t)  = f(u_{0})  & \quad \mathrm{in} \, (0,L)\times(0,\infty),  \\[2mm]		
			v_x(0,t)=0,  v(L,t)=0   & \quad \mathrm{in} \,(0,\infty), \\[2mm]
			v(x,0) =u_0(x)   & \quad \mathrm{on} \, [0,L], 
		\end{cases}
	\end{aligned}
\end{equation}
and its solution $v$, which exists uniquely in $(0, L) \times(0, \infty)$. Similar to the proof of Lemma \ref{lmms0}, 
we have $D_t^{\alpha}v(x,t)>0$  in $(0,L)\times(0,\infty)$. It follows from $u\geq u_0$ and $f(u)\geq f(u_0)$ that 
$u(x,t)>v(x,t)$ in $(0, L) \times(0, T)$.  Let
\begin{equation}\label{dnşmw} \nonumber
w(x,t)=u(x,t)-v(x,t).
\end{equation}
Similar to the proof of Lemma \ref{lmms0}, $w(x,t)$ strictly increases in $t$ and $D_t^{\alpha}w(x,t)\geq 0$. Thus, $D_t^{\alpha}u(x,t)\geq D_t^{\alpha}v(x,t)$ in $(0, L) \times(0, T)$. Let
$$c_1=\min_{[x_0,x_1]\times[t_0,T]}D_t^{\alpha}v(x,t)>0.$$ 
Then $D_t^{\alpha}u(x,t) \geq D_t^{\alpha}v(x,t) \geq c_1>0$ in  $[x_0,x_1]\times\left[t_0,T\right)$.  $\qedsymbol$

The subsequent theorem demonstrates that the quenching points lie in a compact subset of $(0, L)$.

\begin{teo}\label{zk1}
	Suppose that $\Phi_{1}(x)$ is the first eigenfunction of problem (\ref{özd2}).  For any positive constant $\tilde{c}_0<u^{*}$, there exist $x_0$ and $x_1$ with $0<x_0<x_1<L$ such that $u(x,t)<\tilde{c}_0$ in $([0,x_0]\cup[x_1,L])\times \left[0,T\right) $.

\end{teo}
\noindent \textbf{Proof.}  There exists a positive constant $c$ such that
\begin{equation}\nonumber
	\int_{0}^{L} \Phi_{1}(x)(u(x, t))^{2} d x \leq c
\end{equation}
is fulfilled. By considering (\ref{öz6}) and $\int_{0}^{L} \Phi_{1}(x)(u(x,t))^{2} d x < \infty$, we have 
\begin{equation}\nonumber
	\int_{0}^{x} \Phi_{1}(x)dx \rightarrow 0 \quad \text{as} \quad x \rightarrow 0	
\end{equation}
and
\begin{equation}\nonumber
	\int_{0}^{x} \Phi_{1}(x)(u(x,t))^{2} d x \rightarrow 0 \quad \text{as} \quad x\rightarrow 0.	
\end{equation}
Therefore, for any given positive constant $c_0<u^{*}$, there exists a $x_0 \in(0, L/2)$ such that, for $0<x<x_0$,
\begin{equation}\nonumber
	\Bigg(\int_{0}^{x} \Phi_{1}(x) d x\Bigg)^{1 / 2}\Bigg(\int_{0}^{x} \Phi_{1}(x)(u(x,t))^{2} dx\Bigg)^{1 / 2}<c_{0}. 
\end{equation}
We have
\begin{equation}\nonumber
	\begin{aligned}
		\underline{y}(t) &=\int_{0}^{x} \Phi_{1}(x)u(x, t)dx \\[2mm]
		&=\int_{0}^{x} \Phi_{1}(x) \Phi_{1}(x)^{-1 / 2} \Phi_{1}(x)^{1 / 2}u(x, t)dx \\[2mm]
		&\leq\Bigg(\int_{0}^{x} \Phi_{1}(x)dx\Bigg)^{1 / 2}\Bigg(\int_{0}^{x}\Big[\Phi_{1}(x)^{1 / 2} u(x,t)\Big]^{2}dx\Bigg)^{1 / 2} \\[2mm]
		&<c_{0}
	\end{aligned}
\end{equation}
for all $t\in (0,T)$. Similarly, there exists $x_1 \in (L/2,L)$ such that, for $x_1<x<L$, 
\begin{equation}\nonumber
	\Bigg(\int_{x}^{L} \Phi_{1}(x) d x\Bigg)^{1 / 2}\Bigg(\int_{x}^{L} \Phi_{1}(x)(u(x,t))^{2} dx\Bigg)^{1 / 2}<c_{0}. 
\end{equation}
Then we have
\begin{equation}\nonumber
	\begin{aligned}
		\overline{y}(t) &=\int_{x}^{L} \Phi_{1}(x)u(x, t)dx \\[2mm]
		&=\int_{x}^{L} \Phi_{1}(x) \Phi_{1}(x)^{-1 / 2} \Phi_{1}(x)^{1 / 2}u(x, t)dx \\[2mm]
		&\leq \Bigg(\int_{x}^{L} \Phi_{1}(x)dx\Bigg)^{1 / 2}\Bigg(\int_{x}^{L}\Big[\Phi_{1}(x)^{1 / 2} u(x,t)\Big]^{2}dx\Bigg)^{1 / 2} \\[2mm]
		&<c_{0}
	\end{aligned}
\end{equation}	
for all $t\in (0,T)$. Similar to the proof of Theorem \ref{teosnm2}, by denoting $y(t)=\int_{0}^{L} \Phi_{1}(x)u(x, t)dx$, we establish that $y(t)<c_0$ for every $(x,t)\in ([0,x_0]\cup[x_1,L])\times \left[0,T\right)$. 
Hence, considering  (\ref{öz6}) and assumption (5), we conclude that there exists a constant $\tilde{c}_0$ such that $u(x,t) < \tilde{c}_0 $ holds for every $(x,t)\in ([0,x_0]\cup[x_1,L])\times \left[0,T\right)$.
$\qedsymbol$	
	
The following theorem demonstrates the occurrence of the blow-up of $D_t^{\alpha}u(x,t)$ at the quenching point.
	
\begin{teo}
Suppose that $-u_{0}^{''}(x) \leq f(u_0(x))$ in $(0,L)$ and the strict inequality holds somewhere in $(0,L)$. If $\max \{u(x, t): x \in[0, L]\} \rightarrow u^{*}$ as  $t \rightarrow T^{-}<\infty$, where $u^{*}$ is a singular point, then $\sup \{D_t^{\alpha}u(x,t): x \in[0, L]\} \rightarrow \infty$ as $t \rightarrow T^{-}$.
\end{teo}	

\noindent \textbf{Proof.} If $\max \{u(x, t): x \in[0, L]\} \rightarrow u^{*}$ as $t \rightarrow T^{-}<\infty$, then there exist a $x^{*}\in [0,L]$, such that $\lim _{(x, t) \rightarrow\left(x^{*}, T\right)} u(x, t)=u^{*}$. From Theorem \ref{zk1}, we have  $x^{*}\in (x_0,x_1)$. Let 
\begin{equation}\label{dnşms}
	w(x,t)=D_{t}^{\alpha}u(x,t)-\delta f(u) \quad \text{in} \, [x_{0}, x_{1}] \times[0, T)
\end{equation}
where $\delta$ is a positive constant to be determined. Then, we obtain
\begin{eqnarray}
&D_{t}^{\alpha}w(x,t)&=D_{t}^{\alpha}(D_{t}^{\alpha}u(x,t))-\delta D_{t}^{\alpha}f(u) \nonumber \\[2mm]
&\dfrac{\partial}{\partial x}w(x,t)&=\dfrac{\partial}{\partial x}(D_{t}^{\alpha}u(x,t))-\delta f^{\prime}(u) \dfrac{\partial}{\partial x}u(x,t)\nonumber \\[2mm]
&\dfrac{\partial^{2}}{\partial x^2}w(x,t)&=\dfrac{\partial^{2}}{\partial x^2}(D_{t}^{\alpha}u(x,t))-\delta f^{\prime\prime}(u)\Bigg(\dfrac{\partial}{\partial x}u(x,t)\Bigg)^2-\delta f^{\prime}(u) \dfrac{\partial^{2}}{\partial x^2}u(x,t)\nonumber
\end{eqnarray}
Hence,
\begin{eqnarray}
	&D_{t}^{\alpha}w-\dfrac{\partial^{2}w}{\partial x^2}&=D_{t}^{\alpha}\Bigg(D_{t}^{\alpha}u-\dfrac{\partial^{2}u}{\partial x^2}\Bigg)-\delta f^{\prime}(u) D_{t}^{\alpha}u+\delta f^{\prime\prime}(u)\Bigg(\dfrac{\partial u}{\partial x}\Bigg)^2+\delta f^{\prime}(u) \dfrac{\partial^{2}u}{\partial x^2}
	\nonumber \\[2mm]
	&&=D_{t}^{\alpha}f(u)-\delta f^{\prime}(u)\Bigg(D_{t}^{\alpha}u-\dfrac{\partial^{2}u}{\partial x^2}\Bigg)+\delta f^{\prime\prime}(u)\Bigg(\dfrac{\partial u}{\partial x}\Bigg)^2  \nonumber \\[2mm]
	&&=f^{\prime}(u)D_{t}^{\alpha}u-\delta f^{\prime}(u)f(u)+\delta f^{\prime\prime}(u)\Bigg(\dfrac{\partial u}{\partial x}\Bigg)^2 \nonumber \\[2mm]
	&&=f^{\prime}(u)w+\delta f^{\prime\prime}(u)\Bigg(\dfrac{\partial u}{\partial x}\Bigg)^2. \nonumber
\end{eqnarray}
That is,
$$
D_{t}^{\alpha}w-\dfrac{\partial^{2}w}{\partial x^2}-f^{\prime}(u)w= \delta f^{\prime\prime}(u)\Bigg(\dfrac{\partial u}{\partial x}\Bigg)^2\geq 0 \quad \text{in}\, (x_0,x_1)\times(0,T).
$$

Let $\eta\in(0,T)$. Since  $D_{t}^{\alpha}u(x,t)>0$  in $(0,L)\times(0,T)$, it follows from Lemma \ref{lmms1} that there exists a $c_1(>0)$ such that $D_{t}^{\alpha}u(x,t)>c_1$ on $\left\{x_{0}, x_{1}\right\} \times[\eta, T)$, and  $\left[x_{0}, x_{1}\right] \times\{\eta\}$. Theorem \ref{zk1} implies that there exists a $c_{2}<u^{*}$ such that $u(x,t)<c_{2}$ on 
$\left(\left\{x_{0}, x_{1}\right\} \times[\eta, T)\right) \cup\left(\left[x_{0}, x_{1}\right] \times\{\eta\}\right)$. Let us choose 
\begin{equation}\nonumber
	0<\delta<c_{1} / f\left(c_{2}\right).
\end{equation}
Then we can write
$$
w(x, t)=D_{t}^{\alpha}u(x,t)-\delta f(u)>c_{1}-\delta f\left(c_{2}\right)>0 \quad \text{on} \left\{x_{0}, x_{1}\right\} \times[\eta, T) \cup\left[x_{0}, x_{1}\right] \times\{\eta\}.
$$
By the maximum principle, $w(x,t)\geq0$ in $[x_{0}, x_{1}] \times[\eta, T)$. Hence,
$$
\lim _{(x, t) \rightarrow\left(x^{*}, T\right)} D_{t}^{\alpha}u(x, t) \geq \lim _{u \rightarrow u^{{*}^{-}}} \delta f(u)=\infty.
$$
Thus, $\sup \left\{D_{t}^{\alpha}u(x,t): x \in[0, L]\right\} \rightarrow \infty$ as $t \rightarrow T^{-}$ if $\max \{u(x, t): x \in[0, L]\} \rightarrow u^{*}$ as $t \rightarrow T^{-}$. $\qedsymbol$

\begin{uyr}
From the above proof, we can also see that whereever $u(x, t) \rightarrow u^{*}$ as $t \rightarrow T^{-}$, $D_{t}^{\alpha}u \rightarrow \infty$ holds true.
\end{uyr}	
	
\section{Conclusion}
In this study, we prove the occurrence of a blow-up in the time-derivative at the quenching point by employing the Caputo derivative on the time-dependent term.We expect that this analysis will address a significant gap in the literature concerning the fractional quenching problem. Furthermore, we establish these results for a general nonlinear source term, encompassing all potential source terms exhibiting singularities.

\end{document}